\documentclass[dvipdfmx]{article}

\usepackage{amsmath,amsfonts,amssymb,amsthm,pb-diagram,picinpar,graphicx,color}
\usepackage{algorithm}
\usepackage{algorithmic}
\usepackage{colortbl}
\usepackage{enumerate}
\usepackage{hyperref}
\usepackage{bm}
\usepackage{tikz}
\usepackage{pgf}
\usepackage[T1]{fontenc}
\usepackage{enumitem}
\usetikzlibrary{shapes,matrix,calc,arrows,decorations.markings,knots,patterns,intersections,cd}

\newcommand{\CC}{\mathbb{C}}
\newcommand{\PP}{\mathbb{P}}
\newcommand{\mcC}{\mathcal{C}}
\newcommand{\mcB}{\mathcal{B}}
\newcommand{\mcG}{\mathcal{G}}
\newcommand{\bbL}{\mathbb {L}}
\newcommand{\mcO}{\mathcal{O}}

\newcommand{\Kbar}{\overline{K}}

\newcommand{\fd}{\mathfrak{d}}

\newcommand{\pic}{\mathop{\mathrm{Pic}}\nolimits}
\newcommand{\amari}{\mathop{\mathrm{Rem}}\nolimits}
\newcommand{\Div}{\mathop{\mathrm{Div}}\nolimits}
\newcommand{\divi}{\mathop{\mathrm{div}}\nolimits}

\newcommand{\Supp}{\mathop{\mathrm{Supp}}\nolimits}

\newcommand{\sr}{\mathop{\mathrm{sr}}\nolimits}
\newcommand{\red}{\mathop{\mathrm{r}}\nolimits}

\newcommand{\ord}{\mathop{\mathrm{ord}}\nolimits}
\newcommand{\wdeg}{\mathop{\mathrm{wdeg}}\nolimits}
\newcommand{\aff}{\mathop{\mathrm{aff}}\nolimits}

\newcommand{\LT}{\mathop{\mathrm{LT}}\nolimits}
\newcommand{\LM}{\mathop{\mathrm{LM}}\nolimits}
\newcommand{\LC}{\mathop{\mathrm{LC}}\nolimits}
\newcommand{\multi}{\mathop{\mathrm{multideg}}\nolimits}

\newtheorem{thm}{Theorem}[section]
\newtheorem{lem}[thm]{Lemma}     
\newtheorem{cor}[thm]{Corollary}
\newtheorem{prop}[thm]{Proposition}

\theoremstyle{definition}
\newtheorem{defin}[thm]{Definition}

\newtheorem{ex}[thm]{Example}

\theoremstyle{remark}
\newtheorem{rem}[thm]{Remark}

\renewcommand{\thesubparagraph}{\theparagraph.\@arabic\c@subparagraph}
\begin{document}
\begin{center}

{\Large \bf Representations of divisors on hyperelliptic curves, 
 Gr\"obner bases and plane curves with quasi-toric relations}

\bigskip

 {\large \bf Ai Takahashi and Hiro-o Tokunaga\footnote{Partially supported by Grant-in-Aid for Scientific Research C (20K03561)}}
 
 \end{center}

\begin{abstract}
In the study of hyperelliptic curve cryptography, presentations of semi-reduced divisors on a 
hyperelliptic curve play important roles. In this note, we give an interpretation for such presentations
from view points of Gr\"obner bases. As applications,  (i) we give a method to construct weak $n$-contact curves and (ii) we study plane curves satisfying quasi-toric relations of type $(2, n, 2)$
\end{abstract}

\section*{Introduction}

 Let $\mcC$ be a hyperelliptic curve defined
over a field $K$, ${\mathrm{char}}(K) \neq 2$ given by an affine equation
\[
\mcC: y^2 = f(x), \quad f(x) = x^{2g+1} + c_1x^{2g} + \ldots + c_{2g},
\]
where $f(x) = 0$ has no multiple roots in $\Kbar$, where $\Kbar$ is an algebraic closure of $K$. We denote the point at
infinity by $O$.
In the study of hyperelliptic curve cryptography ((\cite{cantor, costello-lauter, galbraith, lauter03}), 
a pair of two polynomials $(u, v)$ $(u, v \in \Kbar[x])$ is used in order to describe semi-reduced divisors  on $\mcC$ (See \S~\ref{subsec:semi-reduced} for
semi-reduced divisors)  and to consider the addition in the Jacobian
of $\mcC$. Such a pair was first considered in \cite{mumfordII} and is
called the Mumford representation of a semi-reduced divisor. 
 For a semi-reduced divisor $\fd$,  $\fd$ is 
given by zeros of  the ideal $\langle u, y -v, y^2  - f\rangle$ generated by $u, y- v, y^2 - f$
in $\Kbar[x, y]$
 with multiplicities.
In \cite{leitenberger}, another description for semi-reduced divisors was given. We  call it the Leitenberger representation.

In this note,  we consider reduced Gr\"obner bases of $\langle u, y - v, y^2 - f\rangle$ with respect to two monomial orders: $1.$ the pure lexicographic order with respect to $y>x$ and $2.$
 a weighted reverse lexicographic order given in
\S~\ref{subsec:monomial-order}. We give interpretations concerning the Mumford representation and the Leitenberger representation
from Gr\"obner bases viewpoints (Propositions~\ref{prop:mumford-rep} and \ref{prop:jacobi-leitenberger}). As applications
we consider the case of $g =1$ and obtain the following results:
\begin{enumerate}

\item[(i)] We give a  method to construct a  weak $n$-contact curve which contacts at a unique $n$-torsion point of $\mcC$ and meets at $O$ (\S~\ref{subsec:n-contact_curve}).
.

\item[(ii)]  Explicit examples of plane curves satisfying infinitely many quasi-toric relations of type $(2, n,2)$
(\S~\ref{subsec:quasi_toric_relations}) are given.

\end{enumerate}

Let us explain these applications briefly. For a smooth cubic $\mcC$, a plane curve $D$ is said to be
a weak $n$-contact curve to $\mcC$ if  the divisor $D|_E$ on $E$ defined by $D$  is of the form
 $D|_E = n \left (\sum_{i=1}^d P_i \right ) +  sO$ for some non-negative integer $s$.
 As we see in \cite{taka-toku20-1}, a weak $n$-contact curve to a cubic as above plays a key role to construct 
examples of  certain Zariski tuples.
As for application (i), we   give
a method to construct weak $n$-contact  curves contact as above, which also works in the case when $n$ is a prime number. In \S~\ref{subsec:examples}, we give explicit examples for $n = 5, 7$.
Note that in \cite{taka-toku20-1}, we only treat with the cases of $n = 3, 4, 6, 8$.

Before we go on to the application (ii), we recall the definition of a quasi-toric relation.
Following to \cite[Definitin 2.13]{jose-anatoly14}, we say that a plane curve $\mcB$ in $\PP^2$ satisfies a
 quasi-toric relation of type $(p, q, r)$ if there exist a sextuple $(F_1, F_2, F_3, h_1, h_2, h_3)$ of
 non-zero homogeneous polynomials such that

\begin{itemize}

\item it satisfies the following relation
\[
h_1^pF_1 + h_2^qF_2 + h_3^rF_3 = 0,
\]
and
\item the curve $\mcB$ is given by $F_1F_2F_3 = 0$.

\end{itemize}
Plane curves that satisfy  quasi-toric relations of certain types has been studied in \cite{jose-anatoly14, kloosterman17, kloosterman18} from the viewpoint of embedded topology of plane curves. 
 We construct examples of curves satisfying infinitely many quasi-toric relations of  type $(2, n, 2)$ $n = 3, 5, 7$ (\S~\ref{subsec:examples}). Note that the cases of $(2, 5,2)$   and $(2, 7, 2)$ were not considered since 
 such cases are not {\it elliptic type} in the terminology of \cite{jose-anatoly14}, i.e., $(2, 3, 6), (3, 3, 3), (2, 4, 4)$.

 \section{Preliminaries}
  
 
\subsection{Two monomial orders on $\Kbar[x, y]$}\label{subsec:monomial-order}

As for general facts on monomial orders and Gr\"obner bases, we refer to \cite{CLO}.
In this note, we consider two monomial orders $>_1$ and $>_2$ as follows:

\begin{itemize}
\item $>_1$ is the pure lexicographic order with $y>x$.

\item $>_2$ is a weighted lexicographic order as follows:
For a monomial $y^mx^n$, we put $\wdeg(y^mx^n) = (2g+1)m + 2n$. We say
$y^{m_1}x^{n_1} >_2 y^{m_2}x^{n_2}$
 if and only if
\begin{enumerate}
\item[(i)]  $\wdeg(y^{m_1}x^{n_1}) > \wdeg(y^{m_2}x^{n_2})$ or 
\item[(ii)] $(2g + 1)m_1 + 2n_1 =  (2g+1)m_2 + 2n_2$ and $n_1 < n_2$
\end{enumerate}
\end{itemize}
 
 The monomial order $>_2$ is nothing but a weighted reverse lexicographic order for $y, x$ with weight $(2g+1, 2)$.
 It coincides with the $C_{ab}$-order considered in \cite{arita} for $(a, b) = (2g+1, 2)$.
  By $\LM_i(g)$,  $\LC_i(g)$ and $\LT_i(g)$, we denote the leading monomial, coefficient and term of $g$ with
  respect to $>_i$, respectively. Also we denote the multidegree with respect to $>_i$ by $\multi_i$.

 \subsection{Semi-reduced divisors on a hyperelliptic curves}\label{subsec:semi-reduced}
 Let $\mcC$ be a hyperelliptic curve defined over $K$ given by the affine equation in the Introduction. We give a summary for
 semi-reduced divisors considered in hyperelliptic cryptocgraphy \cite{costello-lauter, galbraith, lauter03, MWZ} and
 our previous article \cite{taka-toku20-1}. Our notation here are those in \cite{taka-toku20-1}.
 
 Let $\fd$ be a divisor on $\mcC$ and $\Supp(\fd)$ denotes its supporting set. Let $\iota: (x, y) \mapsto (x, -y)$ be
 the hyperelliptic involution on $\mcC$. For any divisor $\fd$ on $\mcC$
  with
$\fd = \sum_{P \in \mcC} m_P P$, 
by considering points of the form $P + \iota(P)$ contained in 
$\fd$, we have a decomposition 
$\fd = \fd_{\sr} + \fd_o$ such that
\begin{enumerate}

\item[(i)] the divisor $\fd_o$ is of the form $\bar\fd + \iota(\bar\fd)$ for some divisor $\bar\fd$, and 

\item[(ii)]  if we write $\fd_{\sr} = \sum_{P \in \mcC}  m'_P P$, then $m'_P$ satisfies the following
conditions:
      \begin{enumerate}
       \item[(a)] $m'_P = 1$ if $m'_P > 0$ and $P = \iota(P)$, and
      \item[(b)] $m'_{\iota(P)} = 0$ if $m'_P > 0$ and $P\neq \iota(P)$. 
       \end{enumerate}
  \end{enumerate}

  We here define a semi-reduced divisor on $\mcC$ following to \cite{galbraith}.
  
  \begin{defin}\label{def:semi-reduced}{\rm Let $\fd$ be a divisor on a hyperelliptic curve $\mcC$ . 
  
  \begin{enumerate}
  
  \item[(i)] The divisor $\fd$ is said to be affine divisor if $\Supp(\fd) \subset \mcC_{\aff}:= \mcC\setminus \{O\}$.
  
  \item[(ii)] An effective affine divisor $\fd$ is said to be  {\it semi-reduced} if $\fd_o$ is empty.  
  \item[(iii)] 
  A semi-reduced divisor $\sum_i m_i P_i$ is said to be $h$-reduced if $\sum_im_i \le g$.
  
  \end{enumerate}
  }
  \end{defin}

\begin{rem}{\rm 
In \cite{galbraith, MWZ}, a semi-reduced divisor satisfying the condition (iii) is simply called a reduced divisor. Here, we use the terminology ^^ {\it $h$-reduced}' in oder to avoid 
confusion for  the terminology {\it reduced divisor} used in standard textbooks in algebraic
geometry e.g., \cite{iitaka}.
}
\end{rem}

Here are some properties for semi-reduced divisors:

\begin{lem}\label{lem:semi-reduced1}{
\begin{enumerate}

 \item[\rm(a)]  For any divisor $\fd = \sum_Pm_PP$ with $\Supp(\fd) \neq \emptyset$, there exists a semi-reduced
divisor $\sr(\fd)$ such that 
 {\rm (i)} $\fd - (\deg \fd)O \sim \sr(\fd) - (\deg\sr(\fd))O$ and
{\rm (ii)} $|\fd| \ge |\sr(\fd)| (= \deg \sr(\fd))$.
Here  we put $\deg \fd := \sum_P m_P$ and
 $|\fd| := \sum_P|m_P|$.

\item[{\rm (b)}] Let $\fd$ be any semi-reduced divisor on $\mcC$ with $\deg \fd > g$. Then there exists a
unique $h$-reduced divisor $\red(\fd)$ such that $\fd - \deg\fd O\sim \red(\fd) - (\deg{\red(\fd)})O$.

\item[{\rm (c)}]  With two statements as above,  we see that for any element $\fd \in \Div^0(\mcC)$, there exists a unique
$h$-reduced divisor $\red(\fd)$ such that $\fd \sim \red(\fd) - (\deg\red(\fd))O$.

\end{enumerate}
}
\end{lem}
As for proofs, see \cite{galbraith, MWZ}.

 \subsection{Representations for semi-reduced divisors}
 
 We keep our notation and terminologies as in \S~\ref{subsec:semi-reduced}. Let 
 $\langle y^2 - f\rangle \subset \Kbar[x,y]$ be the ideal generated by $y^2-f$, where $f$ is the polynomial in 
 the Introduction.
The quotient ring $\Kbar[x, y]/\langle y^2 -f \rangle$ is said to be the coordinate
ring of $\mcC$ and we denote it by $\Kbar[\mcC]$.  The quotient field of $\Kbar[\mcC]$ is the rational function field 
$\Kbar(\mcC)$ of $\mcC$.
 An element of $\Kbar[\mcC]$ is called
a polynomial function, i.e., a rational function with poles only at $O$. 
For $g \in \Kbar[x, y]$, its class in $\Kbar[C]$ gives a
polynomial function on $C$, which we denote by $[g]$. 

For our later use, we define a $\Kbar[x]$-submodule $\amari(y^2)$  of $\Kbar[x, y]$ as follows:
\[
 \amari(y^2)  =  \{b_0(x) + b_1(x)y \mid b_0(x), b_1(x) \in \Kbar[x]\} 
 \]

Since any element in $\Kbar[\mcC]$ can be represented by the class of an element in $\amari(y^2)$ uniquely
(\cite[\S2]{MWZ}), 
we use elements in $\amari(y^2)$ as  normal forms of polynomial functions in $\Kbar[\mcC]$. 
 Let $\fd = \sum_P e_P P$ be a semi-reduced divisor on $\mcC$. We define  ideals $I(\fd) \subset \Kbar[\mcC]$ and $\widetilde{I(\fd)} \subset \Kbar[x, y]$.
\begin{eqnarray*}
I(\fd) &:=& \{[g] \in \Kbar[\mcC] \mid \ord_P([g]) \ge e_P, \mbox{for $\forall P \in \Supp(\fd)$}\}, \\
\widetilde{I(\fd)} & := & \{ g \in \Kbar[x, y] \mid [g] \in I(\fd)\}
\end{eqnarray*}
where $\ord_P$ denotes an order function of the local ring $\mcO_P(\mcC)$ at $P$.  In the remaining of
this section, we consider some generators  of $\widetilde{I(\fd)}$, which are used 
in order to 
 compute the {\it ^^ addition'} on $\pic^0(\mcC)$,
 explicitly, in \cite{galbraith, lauter03, leitenberger, MWZ}.
 
  \subsubsection{Mumford representation}
  
  In \cite{mumfordII}, Mumford gave a description of $\pic^0(\mcC)$ by two polynomials of one variable, by which we have a semi-reduced divisor.
    In the study of hyperelliptic cryptography
  it is called the Mumford representation. We explain it briefly.
  
  Let $\fd = \sum_{i=1}^r e_iP_i (e_i >0)$ be a semi-reduced divisor and put $P_i = (x_i, y_i)$ $(i = 1, \ldots, r)$.

  \begin{lem}\label{lem:mumford-rep}{There exists unique polynomials $u(x), v(x) \in \Kbar[x]$  such that
  \begin{enumerate}
  \item[\rm{(i)}] $u(x) := \prod_{i=1}^r (x - x_i)^{e_i}$, 
  \item[\rm{(ii)}] $\deg v(x) < \deg u(x)$, $\ord_{P_i}([y - v(x)]) \ge e_i$, and  
  \item[\rm{(iii)}] $v(x)^2 - f$ is divisible by $u$.
  \end{enumerate}
  }
  \end{lem}
  For a proof. see \cite[Lemma 10.3.5]{galbraith}.  
  
  \begin{defin}\label{def:mumford-rep}{\rm Let $\fd$ be a non-zero semi-reduced divisor on a hyperelliptic curve $\mcC$. The pair of polynomials
  $(u, v)$ is said to be  the Mumford representation of $\fd$. By $\fd(u, v)$, we mean a non-zero semi-reduced divisor with the Mumford
  representation $(u, v)$. For $\fd = 0$, we take $u(x) = 1$ an $v(x) = 0$ as its Mumford representation.
  }
  \end{defin}
  Note that if $u$ and $v$ as above exist, we recover $\fd$:
  \[
  \fd = \left (\gcd(\divi([u]), \divi([y-v]))\right )_{\aff},
  \]
  where, for $g_i \in \Kbar[x, y]$ $(i = 1, 2)$ and  divisors, $\divi([g_i])$, of functions $[g_i] (i = 1, 2)$, we define
  \[
  \gcd(\divi([g_1]), \divi([g_2])) := \sum_{P\in \mcC_{\aff}} \min(\ord_P([g_1]), \ord_P([g_2]))P -  \left (\sum_P \min(\ord_P([g_1]), \ord_P([g_2])\right) O,
  \]
  and $\left(\gcd(\divi([g_1]), \divi([g_2]))\right )_{\aff} := \sum_{P \in \mcC_{\aff}}\min(\ord_P([g_1]), \ord_P([g_2]))P$.

 As it is shown in \cite{cantor, galbraith, MWZ}, one can compute the addition law on $\pic^0(\mcC)$
 in terms of Mumford representations of two semi-reduced divisors. Also if we are given a semi-reduced
 divisor $\fd(u, v)$,  we have an algorithm to compute the $h$-reduced divisor $\red(\fd(u, v))$ as in Lemma~\ref{lem:semi-reduced1} in terms of $u, v$.

  \subsubsection{Leitenberger representation}
  
  In this subsection, we recall another representation of a non-zero semi-reduced divisor $\fd$ considered in 
  \cite{leitenberger}. 
  Let $\fd = \sum_{i=1}^rP_i$ be
  a non-zero semi-reduced divisor on $\mcC$. By Lemma~\ref{lem:semi-reduced1}, there exists a unique $h$-reduced divisor, $\red(\fd)$,  such that
  \[
  \fd - (\deg\fd) O \sim \red(\fd) - (\deg\red(\fd))O.
  \]
  Hence we have 
  \[
  \fd + \iota^*\red(\fd) - (\deg\fd + \deg \red(\fd))O 
  \sim \red(\fd) + \iota^*\red(\fd)  - 2(\deg\red(\fd)) \sim 0,
  \]
 and there exists a $\psi \in \Kbar[\mcC]$, unique up to constants, such that
 \[
 \divi(\psi) = \fd + \iota^*\red(\fd) - (\deg\fd + \deg \red(\fd))O.
 \]
 Thus we have
 
 \begin{lem}\label{lem:leitenberger0}{ $\deg \red(\fd) = \min\{ r \mid \bbL(-\fd + (\deg\fd + r) O) 
 \neq \{0\}\}$. Here for a divisor $\fd$, $\bbL(\fd)$ denotes vector space consisting of rational functions $\xi$ such that
 $\divi(\xi) + \fd$ is effecitve and $0$.
 }
 \end{lem}

 By choosing $b = b_0 + b_1 y\in \amari(y^2)$ such that $\psi = [b]$, we have

 \begin{lem}\label{lem:leitenberger1}{
 The effective divisor $\fd + \iota^*\!\red(\fd)$ is semi-reduced if and only if
 $\gcd(b_0, b_1) = 1$.
 }
 \end{lem}
 
 \proof 
  Suppose that $\fd + \iota^*\!\red(\fd)$ is not semi-reduced. We then infer 
 that   $\fd + \iota^*\!\red(\fd)$ is of  the form $\fd_1 + P + \iota^*P$ for some effective divisor
 $\fd_1$ and $P=(x_P, y_P)$.  As $P + \iota^*\! P - 2O \sim 0$,
 $\fd_1 - (\deg(\fd + \iota^*\!\red(\fd)) - 2)O \sim 0$. This implies that there exists $\tilde{b} \in \amari(y^2)$ such that 
 $\divi([\tilde{b}]) = \fd_1 - (\deg(\fd + \iota^*\!\red(\fd)) - 2)O$, and we have
  $\divi((x - x_P)\tilde{b}) = \fd + \iota^*\!\red(\fd) - (\deg(\fd + \iota^*\!\red(\fd))O$.
 As $(x - x_P)\tilde{b} \in \amari(y^2)$, $b = c(x - x_P)\tilde{b}$ for some $c \in \Kbar^{\times}$. 
 This means $x - x_P | \gcd(b_0, b_1)$.  Conversely,  if $\gcd(b_0, b_1)$ is not constant, the
 divisor $(\gcd(\divi([b_0]), \divi([b_1)]))_{\aff}$ is contained in $\fd + \iota^*\!\red(\fd)$. Therefore $\fd + \iota^*\!\red(\fd)$ is not semi-reduced.
 \endproof
 \begin{lem}\label{lem:leitenberger2}{Let $\fd = \sum_{i=1}^re_iP_i$ be a semi-reduced divisor such that  
 $\fd + \iota^*\!\red(\fd)$ is semi-reduced. Let $u$ be as in Lemma~\ref{lem:mumford-rep} and let
 $b \in \amari(y^2)$ as above. Then
 $\fd := (\gcd(\divi([u]), \divi([b])))_{\aff}$.
 }
 \end{lem}
 
 \proof  Since $\divi([u]) = \sum_{i=1}^r e_i(P_i + \iota^*\! P_i) - (2\deg\fd) O$, 
 $\divi([b]) =  \fd + \iota^*\red(\fd) - (\deg\fd + \deg \red(\fd))O$ and $\Supp(\fd)\cap \Supp( \red(\fd)) = \emptyset$,
our statement follows.
 
 \endproof
 
 \begin{defin}\label{def:jacobi-leitenberger}{\rm Let $\fd$ be a semi-reduced divisor on $\mcC$ and let $\red(\fd)$
 be the corresponding reduced divisor. Assume that
 \begin{center}
 ($\clubsuit$) $\fd + \iota^*\!\red(\fd)$ is semi-reduced.
 \end{center}
 The pair of polynomials $(u, b), u \in \Kbar[x], b \in \amari(y^2)$ in Lemma~\ref{lem:leitenberger2} is called the Leitenberger representation
 of $\fd$.
 }
 \end{defin}
 
 \begin{rem}{\rm If $\clubsuit$ is not satisfied, i.e., $\fd + \iota^*\!\red(\fd)$ is not semi-reduced, Jacobi's interpolation
 function in \cite{leitenberger} does not seem to give the desired rational function as $\gcd(b_0, b_1)$ in
 Lemma~\ref{lem:leitenberger1} is not $1$.
 }
 \end{rem}

\section{Presentations of semi-reduced divisors and Gr\"obner bases}


We keep our notation and terminology in \S1. Let $\fd$ be a semi-reduced divisor.
The following proposition may be well-known, but we here give its proof.
\begin{prop}\label{prop:mumford-rep} {Let $(u, v)$ be the Mumford representation of $\fd$. Then 
$\widetilde{I(\fd)} = \langle u, y- v \rangle$ and $v^2 - f \in
\langle u \rangle$ in $\Kbar[x]$. In particular, $\{u, y - v\}$ is the reduced Gr\"obner basis of 
$\widetilde{I(\fd)}$ with
respect to  $>_1$.
}
\end{prop}

\proof Since $(u, v)$ is the Mumford representation of $\fd$, by definition, we have
$u, y- v \in \widetilde{I(\fd)}$ and $u | v^2 - f$. In particular, $\langle u, y - v \rangle \subseteq
\widetilde{I(\fd)}$ and $\{y - v, u\}$ is the reduced Gr\"obner basis of $\langle u, y - v\rangle$ by
\cite[Chapter 2]{CLO}. We now show that $\widetilde{I(\fd)} \subseteq \langle u, y-v\rangle$.
Choose any $g \in \widetilde{I(\fd)}$, we apply  \cite[Chapter 2, Thoerem 3 (Division Algorithm)]{CLO} to
our case: $g$ and $F = (y - v, u)$ with respect to  $>_1$. Then we have
\[
g = q_1 (y - v) + q_2 u + r, \quad \mbox{$r \in \Kbar[x], \deg r < \deg u$ if $r \neq 0, q_1, q_2\in \Kbar[x,y]$}.
\]
As $r \in \widetilde{I(\fd)}$, 
\[
\ord_{P_i}([r]) \ge e_{P_i} \quad (\forall P_i = (x_{P_i}, y_{P_i}) \in \Supp(\fd)),
\]
and $r(x_{P_i}) = 0$ for $P_i = (x_{P_i}, 0) \in \Supp(\fd)$. Since $\Kbar[x]$ can be regarded as a subset of $\Kbar[\mcC]$, we infer that $u | r$, i.e., $r = 0$.
Hence $\widetilde{I(\fd)} =  \langle u, y - v\rangle$.

\endproof

\begin{rem}\label{rem:gbasis1}{\rm Our proof of Proposition~\ref{prop:mumford-rep} implies  that any element
 $g \in \Kbar[x]\cap \widetilde{I(\fd)}$ is divisible by $u$, i.e., $\Kbar[x] \cap \widetilde{I(\fd)} =  \langle u \rangle$.
 }
 \end{rem}
 
 We consider  the addition of two semi-reduced divisors $\fd_1$ and $\fd_2$. Assume that
 $\fd_1 + \fd_2$ can be rewritten of the form $\fd_3 + \bar{\fd} + \iota^*\bar{\fd}$, where $\fd_3$ is a non-zero
 semi-reduced divisor and $\bar\fd$ is an effective divisor. Let $(u_3, v_3)$ be the Mumford representation of $\fd_3$ and
 let $u_o \in \Kbar[x]$ be a monic polynomial such that $\divi([u_o]) = \bar\fd + \iota^*\bar\fd - 2(\deg\bar\fd) O$. Then we have
 
 \begin{prop}\label{prop:addition}{Assume that $\fd_3 \neq 0$. Both $u_ou_3$ and $u_o(y - v_3)$ are contained in the reduced Gr\"obner basis of 
 $\widetilde{I(\fd_1+\fd_2)}$ with respect to $>_1$.
 }
 \end{prop}
 
 \proof If $\bar\fd = 0$, we can take $1$ as $u_o$. Hence our statement follows from Proposition~\ref{prop:mumford-rep}. Now we assume $\bar\fd \neq 0$.  Put $I_3 = \widetilde{I(\fd_1+\fd_2)}$.  Since $\Kbar[\mcC]$ is a Dedekind domain, we have $I(\bar\fd + \iota^*\bar\fd) = \langle [u_o]\rangle$ and 
 \[
 I(\fd_1 + \fd_2) = I(\fd_3)I(\fd_o + \iota^*\fd_o) = \langle [u_o][u_3], [u_o][y - v_3]\rangle.
 \]
 Hence $I_3  = \langle u_ou_3, u_o(y - v_3), y^2 -f \rangle$. Since
 $\fd_3,  \bar\fd \neq 0$, for any element $g$ in $I_3$,  $\divi([g]) - \bar\fd - \iota^*\bar\fd$ is effective.
 Therefore  no polynomial of the form $y - b, b \in \Kbar[x]$, is contained in $I_3$.
 As $y^2 - f \in I_3$, by \cite[Chapter 5, \S3]{CLO}, 
 the reduced Gr\"obner basis $\mcG(I_3)$ of $I_3$ is of the form $\{g_1, g_2, g_3\}$ such that
 $\LT_1(g_1) = x^{n_1}, \LT_2(g_2) = x^{n_2}y, \LT_1(g_3) = y^2$.  
 We first show that  $g_1 = u_ou_3$, i.e., $I_3 \cap \Kbar[x] = \langle u_ou_3 \rangle$. As $g_1 \in I_3$
 and $\divi([g_1]) - \divi([u_o][u_3])$ is effective,
 $g_1 = u_ou_3g'_1 + h(y^2 -f) $ for some $g'_1 \in \amari(y^2), h \in \Kbar[x, y]$. 
  Since $g_1 \in \Kbar[x]$,  we infer that $g'_1 \in \Kbar[x]$ and $h = 0$. On the other hand, as  $u_ou_3 \in I_3 \cap \Kbar[x]$,  $x^{n_1}$ divides $\LT_1(u_ou_3)$. Hence we have $\LT_1(g_1) = \LT_1(u_ou_3)$ and  $g_1 = u_ou_3$.  
 We next consider $g_2$. As $g_2 \in I_3$,  $g_2 - u_og'_2 \in \langle y^2 -f \rangle$ for some $g'_2 \in \amari(y^2)\setminus \Kbar[x]$. This means that $y^2$ divides $\LT_1(g_2 - u_og'_2)$. 
 As $g_2 \in \amari(y^2)$, we infer that $g_2 = u_og'_2$. This implies $\LT_1(u_o(y - v_3)) = x^{\deg u_o}y$ divides $\LT_1(g_2) = x^{n_2}y$.
 On the other hand, since  $u_o(y - v_3) \in I_3$,  $x^{n_2}y$ divides $\LT_1(u_o(y - v_3))$. This implies that $\LT_1(u_o(y - v_3))
 = \LT_1(g_2)$ and we have
 $u_o(y - v_3) - g_2\in I_3\cap \Kbar[x]$. Hence we have $u_o(y - v_3) = g_2 + rg_1 $ for some $r \in \Kbar[x]$. As $g_1, g_2 \in
 \mcG(I_3)$ and $\deg v_3 < \deg u_3$,    we infer that $\deg rg_1 < \deg g_1$ if $r \neq 0$. Thus  we have $r = 0$ and
 $g_2 = u_o(y - v_3)$.
 \endproof
 
 By our proof of Proposition~\ref{prop:addition},  we have the following corollary.
 \begin{cor}\label{cor:gbasis}{If we let $f \equiv f_o \mod u_ou_3$, then 
 $\{u_ou_3, u_o(y - v_3), y^2 - f_o\}$ is the reduced Gr\"obner basis of $\widetilde{I(\fd_1 + \fd_2)}$
 with respect to $>_1$.}
 \end{cor}
 
 \begin{rem}\label{rem:gbasis1}{\rm 
 By Proposition~\ref{prop:addition}, if  $\fd_1$ and $\fd_2$ are given by the Mumford
 representations $(u_1, v_1)$ and $(u_2, v_2)$, respectively, we obtain the Mumford representation of $\fd_3$ by
 computing the reduced Gr\"obner basis of $\langle u_1u_2, u_1(y - v_2), u_2(y - v_1), (y - v_1)(y - v_2), y^2 - f\rangle$ with respect to $>_1$.
 }
 \end{rem}
  
 We next consider the case for the monomial order $>_2$.
 
 \begin{lem}\label{lem:gbasis-L1}{Let $\fd$ be a semi-reduced divisor and $\red(\fd)$ denotes the unique reduced divisor as in \ref{subsec:semi-reduced}.  Put $w_o = \min\{\wdeg(\LM_2(g)) \mid g \in \widetilde{I(\fd)}, [g] \neq 0\}$. Then
 \[
 \deg\fd + \deg\red(\fd) = w_o
  \] 
 holds.
 }
 \end{lem}
 
 \proof Our proof consists of $2$ steps.
 
 \underline{Step 1}.  We show that $w_{\deg}:=\min\{\multi_2(g) \mid g \in \widetilde{I(\fd)}, [g] \neq 0\}$ is attained by some
 elements in $\widetilde{I(\fd)}\cap \amari(y^2)$.
In particular, $w_{o} = \min\{\wdeg(\LM_2(b)) \mid b \in \widetilde{I(\fd)}\cap \amari(y^2), [b] \neq 0\}$
 
  Choose $g \in \widetilde{I(\fd)}$ arbitrary. Note that $g$ can be expressed uniquely as follows:
 \[
 g = q_g(y^2 - f) + b_g,\,\,  q_g \in \Kbar[x, y], \, b_g= b_{0, g} + b_{1, g}y \in \amari(y^2).
 \]
 As $\multi_3(q_g(y^2 -f)) = \multi_2(g_g) + (2, 0)$ and $\multi_2(b_g) = \max\{(0, \deg b_{0,g}), (1, \deg b_{1, g})\}$, 
 by \cite[Lemma 8, Ch.2]{CLO}, 
 \begin{eqnarray*}
 \multi_2(g) & = & \max\{\multi_2(q_g(y^2 -f)), \multi_2(b_g)\} \\
                  & \ge & \multi_2(b_g)
  \end{eqnarray*}
  Hence $w_{\deg}$ is attained by some element $b$ in $\widetilde{I(\fd)}\cap \amari(y^2)$ with $[b] \neq 0$.
  
  \underline{Step 2}. We first recall that, for $b = b_0 + b_1y \in \amari(y^2)$, $\wdeg(\LT_3(b)) = \max\{2\deg b_0, 2\deg b_1 + (2g+1)\}$ and 
  $\wdeg(\LT_3(b)) = -\ord_O([b])$ holds. Choose $b_{\min}:= b_{0, \min} + b_{1, \min} y$ such that
  $w_{\deg} = \multi_2 b_{\min}$. Then by the definition of $>_2$, we have $w_o = \wdeg(\LT_2(b_{\min}))$.
  Let $b_{\fd}$ be the rational function as in Lemma~\ref{lem:leitenberger1}.
  Since $\fd$ is effective, by Lemmas~\ref{lem:semi-reduced1} and \ref{lem:leitenberger0},  $[b_{\fd}] = [cb_{\min}]$ for some $c \in K^{\times}$ and $\ord_O b_{\fd} = \ord_O b_{\min}$. Hence our statement follows.
 \endproof

 Choose $b_{\fd} = b_0 + b_1y$ as above such that $\LC_2(b_{\fd}) = 1$.
  Let $\mcG_2$ be the reduced Gr\"obner basis of $\widetilde{I(\fd)}$ with
respect to  $>_2$.

\begin{lem}\label{lem:gbasis-L2}{$b_{\fd}$ is a member of  $\mcG_2$.}
\end{lem}

\proof Choose $g_o \in \mcG_2$ so that $\multi_2(g_o)$ is minimum among polynomials in $\mcG_2$ with
$[g_o] \neq 0$.  By our choice of $b_{\fd}$, 
$\multi_2(b_{\fd}) \le \multi_2(g_o)$. On the other hand,  $\LT_2(b_{\fd})$ is divisible by $\LT_2(g)$ for some $g \in \mcG_2$ with $[g] \neq 0$ as $b_{\fd} \in \amari(y^2)$. 
Hence $\multi_2(b_{\fd}) \ge \multi_2(g_o)$ and this implies $\multi_2(b_{\fd}) = \multi_2(g_o)$. Since $b_{\fd}, g_o \in \amari(y^2)$ and both leading coefficients of $b_{\fd}$ and $g_o$ is $1$,
if $b_{\fd} \neq g_o$, then we have $b_{\fd} - g_o \in \widetilde{I(\fd)}$, $[b_{\fd} - g_o ] \neq 0$ and   
$\multi_2(b -g_o) < \multi(b_{\fd})$. This contradicts to our choice of $b_{\fd}$.
Hence $b_{\fd} = g_o \in \mcG_2$.
\endproof

Now we have the following proposition, 

\begin{prop}\label{prop:jacobi-leitenberger}{Let $\fd$ be a semi-reduced divisor satisfying $\clubsuit$. Let $(u, v)$ be  the Mumford
representation of $\fd$ and let $b_{\fd}$ be as above. Then we have $\widetilde{I(\fd)} = \langle u, b_{\fd}, y^2 -f \rangle$. In 
particular, $I(\fd) = \langle [u], [b_{\fd}]\rangle$}
\end{prop}

\proof  Put $b_{\fd} = b_0 + b_1 y$. As $\fd$ satisfies $\clubsuit$, $\gcd(b_0, b_1) = 1$. If $\gcd (u, b_1) \neq 1$, there exists $P = (x_P, y_P) \in \mcC$
such that $u(x_P) = 0$ and $b_1(x_P) = 0$. As $b_{\fd} \in \widetilde{I(\fd)}$, we infer that $b_{\fd}(x_P, y_P) = 0$, i.e., $b_0(x_P) = 0$ also holds. This contradicts to $\gcd(b_0, b_1) = 1$. Hence
 $\gcd (u, b_1) =1$. By choosing  $h_1, h_2 \in \Kbar[x]$ such that $h_1u + h_2b_1 = 1$, we have
 \[
 h_1uy + h_2b_{\fd} = y + h_2b_0.
 \]
 Take $v_1$ so that $-v_1 \equiv h_2b_1 (\bmod u)$.   As 
 $y - v, y- v_1 \in \widetilde{I(\fd)}$, we infer that $v_1 - v \in \widetilde{I(\fd)}\cap \Kbar[x]$. By Remark~\ref{rem:gbasis1}, 
 $v_1 - v$ is divisible by $u$, which implies $v = v_1$. Hence $\langle u, y - v \rangle \subseteq \langle u, b_{\fd}, y^2 -f \rangle$ and 
 our statement follows from Proposition~\ref{prop:mumford-rep}.
\endproof

\section{The case of $g=1$}


In this section, we consider the case of $g = 1$ and apply our results on Leitenberger representations 
in \S 2 to study explicit construction
of plane curves. In this case, $\mcC$ is an elliptic curve and we denote it by
\[
E: y^2 = f(x) = x^3 + ax^2 + bx + c, \quad a, \, b, \, c \in K.
\]

 \subsection{Weak $n$-contact curves to $E$}\label{subsec:n-contact_curve}

Let   $T= (x_T, y_T) \in E(K)$ be a torsion point of order $n$. 
In our previous article \cite{taka-toku20-1}, an element $\xi \in \Kbar[E]$ such
that $\divi(\xi) = n(T - O)$ plays an important part. We give a method in construction $b_{nT} \in \amari(y^2)$ such that
$[b_{nT}] = \xi$ explicitly. 
Let $[k]T$ be the multiplication-by-$k$ of $T$ on $E$. Since $T$ is order $n$, $[k]T \neq T$ for $2 \le k \le n$. This means
$kT$ is a semi-reduced divisor on $E$ satisfying $\clubsuit$. As the Mumford representation of $T$ is $(x - x_T, y_T)$,
${I(T)} = \langle [x - x_T],[ y - y_T]\rangle$ and we infer that $\widetilde{I(kT)} = \langle (x - x_T)^k, (x - x_T)^{k-1}(y - y_T),
\ldots, (y- y_T)^k, y^2 - f \rangle$ for $1 \le k \le n$. 
Compute the reduced Gr\"obner basis $\mcG_2$ of $\widetilde{I(kT)}$ with respect to $>_2$. 
 Let $g_k$ be an element of $\mcG_2\setminus \{y^2  -f\}$ with minimum multidegree. By our proof of Lemma~\ref{lem:gbasis-L2}, $g_k \in \amari(y^2)$ and
 $\divi([g_k]) = kT +[-k]T - (k+1)O$. 
 
 We now consider the case of $k = n$. As $[n]T = O$,  $\divi([g_n]) = nT  - nO$. Hence we can choose $g_n$ as $b_{nT}$. Now put
 $b_{nT} = b_0 + b_1 y$. We have
 \[
(\ast) \quad  b_{nT}\iota^*b_{nT} = b_0^2 - b_1^2y^2 = b_0^2 - b_1^2f = r(x - x_T)^n, \quad r \in K^{\times}
 \]
 as $\divi([ b_{nT}\iota^*b_{nT}]) = nT + n[-1]T - 2nO = \divi([(x - x_T)^n])$. 
 
 Now assume that $K = \CC$ and  let $b_{nT}$ be as above. As $\divi([b_{nT}]) = nT - nO$, the plane curve
 given by $b_{nT} = 0$ is a weak $n$-contact curve to $E$ whose contact affine point is $T$ only.

  \subsection{Curves with quasi-toric relation of  type $(2, n, 2)$}\label{subsec:quasi_toric_relations}
  
  We next consider curves with many quasi-toric relations.
 To this purpose, we generalize the  observation in \S~\ref{subsec:n-contact_curve}
 to a semi-reduced divisor $\fd_T$ such that 
 (i) $\fd_T - \deg(\fd_T)O
 \sim T - O$ and (ii)  $[k]T \not\in \Supp(\fd_T)$ for $1\le k \le n-1$. Since $\red(k\fd_T) = [k]T$,
 $k\fd_T + \iota^*\red(k\fd_T)$  is semi-reduced, i.e., $k\fd_T$ satisfies the condition $\clubsuit$. 
 Let $(u_{k\fd_T}, v_{k\fd_T})$ be the Mumford
 representation of $\fd_T$.  We have ${I(T)} = \langle [u_{\fd_T}],[ y - v_{\fd_T}]\rangle$  and 
 $I(k\fd_T) = I(\fd_T)^k
  = \langle [u_{\fd_T}]^k, [u_{\fd_T}]^{k-1}[(y - v_{\fd_T})], \ldots [(y - v_{\fd_T}]^k\rangle$. 
  Hence 
  $\widetilde{I(k\fd_T)} =  \langle u_{\fd_T}^k, u_{\fd_T}^{k-1}(y - v_{\fd_T}), \ldots, 
  (y - v_{\fd_T})^k, y^2 - f\rangle$. 
  We again compute the reduced Gr\"obner basis $\mcG_2$ of $\widetilde{I(kT)}$ with respect to $>_2$. 
 Let $g_k$ be an element of $\mcG_2\setminus \{y^2  -f\}$ with minimum multidegree. By our proof of Lemma~\ref{lem:gbasis-L2}, 
 $g_k \in \amari(y^2)$ and
 $\divi([g_k]) = k\fd_T +[-k]T - (k\deg{\fd_T}+1)O$. For the case of $k = n$, as $[-n]T = O$, we have
 $\divi([g_n]) = n\fd_T  - (n\deg\fd_T)O$ and we can choose $g_n$ as $b_{n\fd_T}$. 
 Now put
 $b_{n\fd_T} = b_0 + b_1 y$. We have
 \[
(\ast) \quad  b_{n\fd_T}\iota^*b_{n\fd_T} = b_0^2 - b_1^2y^2 = b_0^2 - b_1^2f = r(x - u_{\fd_T})^n, \quad r \in K^{\times}
 \]
 as $\divi([ b_{n\fd_T}\iota^*b_{n\fd_T}]) = n\fd_T + n\iota^*\fd_T - 2n(\deg\fd_T)O = \divi([(x - u_{\fd_T})^n])$. 
   Now assume that $K = \CC(t)$ and $a, b, c \in \CC[t]$. 
  In this case, $f(x) \in \CC[t,x]$ and  we have a plane curve $\mcB_o$ in $\PP^2$ given by the affine equation $f(x) = 0$.
  Put  $\mcB = \mcB_o + L_{\infty}$, where $L_{\infty}$ denotes the line at infinity. 
  If  $r \in \CC^{\times}$ and $u_{n\fd_T}, b_{n\fd_T}$ are also in $ \CC[t, x]$, by homogenizing  $(\ast)$, we see that
  either $\mcB_o$ or $\mcB$ satisfies  a quasi-toric relation of type $(2, n, 2)$.  
Based on this approach, we construct examples of plane curves satisfying  infinitely many quasi-toric relations 
of type $(2, n, 2)$ in \S~\ref{subsec:examples}
 \begin{rem}{By \cite[Theorem 5.1]{cox-parry},  our method using the cases of genus $1$ works 
 for $n \le 12$. 
 }
 \end{rem}
   
 \subsection{Examples}\label{subsec:examples}

 We consider examples for $n = 3, 5, 7$. Note that the latter two  cases are not considered in \cite{taka-toku20-1}.
 
 \begin{ex}\label{eg:3tor}{
 Let $E_3$ be an elliptic curve defined over $K$ given by
$y ^2 = f_3, f_3(x) = x^3 + (mx+ n)^2, m, n \in \CC[t]$. $T:= (0, n)$ is a point  on $E_3$ of oder $3$. 
Note that the tangent line at $T$ is a weak $3$-contact curve to $E$ for general $m, n \in \CC$.
Let $l_{[-1]T}$ be a line through $[-1]T$. We may assume that $l_{[-1]T}$ is given by
$y = rx -n$. Put $l_{[-1]T}\cap E_3 =\{[-1]T, P_1, P_2\}$ and $\fd_T:= P_1 + P_2, P_i = (x_i, y_i)$.
The Mumford representation of $\fd_T$ is $(u_{\fd_T}, v_{\fd_T})$, where
\[
u_{\fd_T} = x^2 + (m^2 - r^2)x + 2n(m-r), \quad v_{\fd_T} = rx -n.
\]
For a general $r$, $\fd_T$ satisfies the condition $\clubsuit$. Now we apply
our argument above to $\widetilde{I(3\fd)}$. Then we have $b_{3\fd_T} = b_0 + b_1y$
where
\begin{eqnarray*}
b_0 & = & x^3 + (2m^2 + 3mr + 3r^2)x^2 + (m^4 + 3m^3r + 3m^2r^2 + mr^3 + mn - 3nr)x \\
& & + m^3n + 3m^2nr + 3mnr^2 + nr^3 + 2n^2 \\
b_1 & = & x(-m - 3r) - m^3 - 3rm^2 - 3mr^2 - r^3 + 2n \\
\end{eqnarray*}
and
\[
b_0^2 - b_1^2f _3= u_{\fd_T}^3.
\]
Now assume that $K = \CC(t)$ and  choose $m, n \in \CC[t]$. Then the above relation
can be considered as (affine ) quasi-toric relation of type $(2, 3, 2)$. In fact,
put $m = 1, n = t, r\in \CC$. Then we have
\begin{eqnarray*}
b_0 & = & x^3 + (2 + 3r + 3r^2)x^2 + (1 + r + 3r^2 + r^3 + t - 3tr)x \\
& & + t + 3rt + 3r^2t + r^3t + 2t^2 \\
b_1 & = & x(-1 - 3r) - 1 - 3r - 3r^2 - r^3 + 2t 
\end{eqnarray*}
and
\[
b_0^2 - b_1^2(x^3 + (x + t)^2) = (x^2 + (1-r^2)x + 2(1-r)t)^3.
\]
By homogenizing both hand side $[T, X, Z], t = T/Z, x = X/Z$, we have
\[
(Z^3b_0(T/Z, X/Z))^2 + (Z(b_1(T/Z, X/Z))^2(Z(X^3 + (X+ tZ)^2Z) = (X^2 + (1-r^2)XZ + 2(1-r)TZ)^3.
\]
Since we can choose $r$ arbitrary, $\mcB$ given by $Z(X^3 + (X+ T)^2Z)= 0$ satisfies
infinitely many quasi-toric relations of type $(2, 3, 2)$ such that 
$F_1 = 1, F_2 = -1, F_3 = Z(X^3 + (X+ T)^2Z), h_1 = Z^3b_0(T/Z, X/Z), h_2 = X^2 + (1-r^2)XZ + 2(1-r)TZ$,
and $h_3 = Z(b_1(T/Z, X/Z)$.

 }
 \end{ex}
 \begin{ex}\label{eg:5tor}{Let $E_5$ be an elliptic curve over $\CC(t)$ given by
 \[
E_5: y^2 = f_5(t,x) = x^3 + \frac14 (t^2 + 4t -4)x^2 + \frac12 t(t - 1)x + \frac 14(t - 1)^2.
\]
Put $\displaystyle{T:=  \left [0, \frac{t-1} 2\right ]}$.  $T$ is a point on $E_5$ of order $5$.
 Then $\widetilde{I(5T)} = \langle x^5, x^{4}(y -  \frac{t-1} 2),  x^{3}(y -  \frac{t-1} 2)^2,  x^{2}(y -  \frac{t-1} 2)^3,
  x(y -  \frac{t-1} 2)^4, (y -  \frac{t-1} 2)^5, y^2 - f_5\rangle$. We have the reduced Gr\"obner basis $\mcG_2(\widetilde{I(5T)})$ of $\widetilde{I(5T)}$ with
  respect to $>_2$ and $b_{5T}$ is as follows:
  \[
  \mcG_2(\widetilde{I(5T)})=\{g_1, g_2, g_3\},
  \]
  where
  \begin{eqnarray*}
  g_1 &= & (-t - 2)x^2 + 2xy + (-2t + 1)x + 2y + 1 - t, \\
  g_2 &= &-4x^3 + (-t^2 - 4t + 4)x^2 + (-2t^2 + 2t)x + 4y^2 - t^2 + 2t - 1,\\
  g_3 & = & 2x^4 - 2x^3 + tx + 2x^2 + t - 2y - 1
  \end{eqnarray*}
  Since $g_2 = 4(y ^2 - f_5)$,  $\multi_2(g_1) = (1,1), \wdeg(\LM_2(g_1)) = 5$ and
  $\multi_2(g_3) = (4, 0), \wdeg(\LM_2(g_3)) = 8$, we have
  $b_{5T} = g_1$. Hence $b_{5T}:=b_0 +b_1 y, b_0:=(-t-2) x^2+(-2 t+1) x+1-t, b_1:= 2x + 2$  and we have
  \[
  b_0^2 - b_1^2f_5 = - 4x^5.
  \]
   For a general $t \in \CC$, the curve $D_{5T}$ given by $b_{5T} = 0$ is a weak $5$-contact curve to $E_5$
   such that $D_{5T}|_{E_5} = 5T + O$. In particular, $\divi([b_{5T}]) = 5T - 5O$.
  
  We next consider curves with quasi-toric relations of type $(2, 5, 2)$.
  We choose  any semi-reduced divisor $\fd_T$ of degree $2$ such that $v_{\fd_T}$ in the
  Mumford representation $(u_{\fd_T}, v_{\fd_T})$ is of the form $r(x - x_T) - y_T, \, r \in \CC^{\times}$. We infer that
  $u_{\fd_T}$ and $v_{\fd_T}$ satisfy 
  \[
  f - v_{\fd_T}^2 = (x - x_T)u_{\fd_T}
  \]
  \[
u_{\fd_T} = -x^2 + (r^2 - \frac{1}{4}t^2 - t + 1)x - rt + r - \frac{t^2}{2} + \frac{t}{2}, \quad v_{\fd_T} = rx-\frac{1}{2}(t-1).
\]
For a general $r$, $\fd_T$ satisfies the condition $\clubsuit$. Now we apply
our argument  to $\widetilde{I(5\fd)}$. Then we have $b_{5\fd_T} = b_0 + b_1y$
where
\begin{eqnarray*}
b_0 & = & 32 - 128t + 32r^5t + 80r^4t^2 + 80r^3t^3 + 40r^2t^4 + 10rt^5 - 80r^4t - 80r^3t^2 - 40r^2t^3 - 10rt^4\\
&& - 128t^3 + 32t^4 - 32r^5 - t^5 + 64x^5 + t^6 + (640r^2 + 160rt + 48t^2 + 320r + 192t - 128)x^4\\
&&  + (320r^4 + 320r^3t + 320r^2t^2 + 80rt^3 + 12t^4 + 640r^3 + 1280r^2t + 480rt^2 + 96t^3 - 640r^2\\
&& + 160t^2 - 160r - 368t + 160)x^3 + (32r^5t + 80r^4t^2 + 80r^3t^3 + 40r^2t^4 + 10rt^5 + t^6 + 64r^5\\
&& + 320r^4t + 480r^3t^2 + 320r^2t^3 + 100rt^4 + 12t^5 + 320r^3t + 800r^2t^2 + 240rt^3 + 56t^4 - 320r^3\\
&& - 1120r^2t - 400rt^2 + 8t^3 + 320r^2 + 320rt - 208t^2 - 160r + 208t - 64)x^2 + (64r^5t + 160r^4t^2\\
&& + 160r^3t^3 + 80r^2t^4 + 20rt^5 + 2t^6 - 32r^5 + 80r^4t + 240r^3t^2 + 200r^2t^3 + 70rt^4 + 9t^5 - 160r^4\\
&& - 320r^3t - 240r^2t^2 - 240rt^3 + 6t^4 + 480rt^2 - 16t^3 - 480rt - 48t^2 + 160r + 80t - 32)x + 192t^2\\
b_1 & = & (-320r - 32t - 64)x^3 + (-640r^3 - 320r^2t - 160rt^2 - 16t^3 - 640r^2 - 640rt - 96t^2 + 320r\\
&& + 32t)x^2 + (-64r^5 - 160r^4t - 160r^3t^2 - 80r^2t^3 - 20rt^4 - 2t^5 - 320r^4 - 640r^3t - 480r^2t^2\\
&& - 160rt^3 - 20t^4 - 320rt^2 - 32t^3 + 640rt + 128t^2 - 320r - 160t + 64)x - 64r^5 - 160r^4t\\
&& - 160r^3t^2 - 80r^2t^3 - 20rt^4 - 2t^5 + 64t^3 - 192t^2 + 192t - 64 \\
\end{eqnarray*}
and
\[
b_0^2 - b_1^2f _5= -4u_{\fd_T}^5.
\]
By homogenizing both hand side $[T, X, Z], t = T/Z, x = X/Z$, we have
\begin{eqnarray*}
&&(Z^8b_0(T/Z, X/Z))^2 + (Z^6b_1(T/Z, X/Z))^2(Z^4f_5(T/Z, X/Z))\\
&& = 4Z(-4XZ^2r^2 + 4TZ^2r - 4Z^3r + T^2X + 2T^2Z + 4TXZ - 2TZ^2 + 4X^2Z - 4XZ^2)^5.
\end{eqnarray*}
Since we can choose $r$ arbitrary, $\mcB$ given by $-4Z^5(f_5(T/Z, X/Z))= 0$ satsifies
infinitely many quasi-toric relations of type (2, 5, 2) such that 
\begin{eqnarray*}
h_1&=&Z^8b_0(T/Z, X/Z),\\
h_2&=&-4XZ^2r^2 + 4TZ^2r - 4Z^3r + T^2X + 2T^2Z+ 4TXZ - 2TZ^2 + 4X^2Z - 4XZ^2, \\
h_3&=&Z^6b_1(T/Z, X/Z), F_1=1, F_2=-4Z \,  \mbox{and}\,   F_3=Z^4f_5(T/Z, X/Z).
\end{eqnarray*}
 }
 \end{ex}
 
 \begin{ex}\label{eg:7tor}{Let $E_7$ be an elliptic curve over $\CC(t)$ given by
 \[
E_7: y^2 = f_7(t, x) = x^3 + \frac 14 (t^4 - 6t^3 + 3t^2 + 2t +1)x^2 + \frac 12 (t^5 - 2t^4 +t^2)x + \frac 14(t^6 - 2t^5 + t^4).
\]
Put $\displaystyle{T:=  \left [0, \frac{t^3- t^2} 2\right ]}$.  $T$ is a point on $E_5$ of order $7$. Then
$\widetilde{I(7T)} = \langle x^7, x^6(y -  \frac{t^3- t^2}2 ),  x^5(y -   \frac{t^3- t^2}2)^2,  x^4(y -   \frac{t^3- t^2}2)^3,
  x^3(y -   \frac{t^3- t^2}2)^4, x^2(y -   \frac{t^3- t^2}2)^5,  x(y -   \frac{t^3- t^2}2)^6,  (y -   \frac{t^3- t^2}2)^7, y^2 - f_7\rangle$. 
  We now compute the reduced Gr\"obner basis $\mcG_2(\widetilde{I(7T)})$  of $\widetilde{I(7T)}$ with
  respect to $>_2$ and $b_{7T}$ is as follows:
  \[
  \mcG_2(\widetilde{I(7T)})=\{g_1, g_2, g_3\},
  \]
  where
  \begin{eqnarray*}
  g_1 &= & -4x^3 + (-t^4 + 6t^3 - 3t^2 - 2t - 1)x^2 + (-2t^5 + 4t^4 - 2t^2)x + 4y^2 - t^6 + 2t^5 - t^4,\\
  g_2 &= &(-t^2 + 3t + 3)x^3 + 2x^2y + (-3t^3 + 4t^2 + 3t + 1)x^2 + (4t + 2)yx + (-3t^4 + 2t^3 + 2t^2)x \\
  &&+ 2yt^2 - t^5 + t^4,\\
  g_3 & = & 2x^5 + (-4t - 2)x^4 + (6t^2 + 8t + 2)x^3 + (-t^5 - 3t^4 + 8t^3 + 12t^2 + 6t+ 1)x^2 + (2t^3+ 12t^2\\
  && + 10t + 2)yx + (-2t^6 - 7t^5 + 4t^4 + 8t^3 + 2t^2)x + (2t^4 + 6t^3 + 2t^2)y- t^7 - 2t^6 + 2t^5 + t^4
  \end{eqnarray*}
  
  Since $g_1 = 4(y ^2 - f_7)$,  $\multi_2(g_2) = (2,1), \wdeg(\LM_2(g_2)) = 7$ and
  $\multi_2(g_3) = (5, 0), \wdeg(\LM_2(g_3)) = 10$, we have
  $b_{7T} = g_2$. Hence $b_{7T}:=b_0 +b_1 y, b_0:=(-t^2 + 3t + 3)x^3 + (-3t^3 + 4t^2 + 3t + 1)x^2 + (-3t^4 + 2t^3 + 2t^2)x - t^5 + t^4, b_1:= 2x^2 + (4t + 2)x + 2t^2$  and we have
  \[
  b_0^2 - b_1^2f_7 = - 4x^7.
  \]
 For a general $t \in \CC$, the curve  $D_{7T}$ given by $b_{7T} = 0$ is a weak $7$-contact curve to $E_7$ such that $D_{7T}|_{E_7} = 7T + 2O$. In particular, $\divi([b_{7T}]) = 7T - 7O$.
 
 We next consider curves with quasi-toric relations of type $(2, 7, 2)$.
 We first choose  any semi-reduced divisor $\fd_T$ of degree $2$ such that $v_{\fd_T}$ in the
  Mumford representation $(u_{\fd_T}, v_{\fd_T})$ is of the form $r(x - x_T) - y_T, \, r \in \CC^{\times}$. We infer that
  $u_{\fd_T}$ and $v_{\fd_T}$ satisfy 
  \[
  f - v_{\fd_T}^2 = (x - x_T)u_{\fd_T}
  \]
   \[
u_{\fd_T} = -x^2 + (r^2 - \frac{1}{4}t^4 + \frac{3}{2}t^3 - \frac{3}{4}t^2 - \frac{1}{2}t - \frac{1}{4})x - rt^3 + rt^2 - \frac{t^5}{2} + t^4 - \frac{t^2}{2},
\quad v_{\fd_T} = rx - \frac{1}{2}t^3 + \frac{1}{2}t^2.
\]
For a general $r$, $\fd_T$ satisfies the condition $\clubsuit$. Now we apply
our argument to $\widetilde{I(7\fd)}$. Then we have $b_{7\fd_T} = b_0 + b_1y$ and
\[
b_0^2 - b_1^2f _7= -4u_{\fd_T}^7.
\]
We here omit explicit forms of $b_0$ and $b_1$ as they are too long.
By homogenizing both hand side $[T, X, Z], t = T/Z, x = X/Z$, we have
\begin{eqnarray*}
&&(Z^{19}b_0(T/Z, X/Z))^2 + (Z^{16}b_1(T/Z, X/Z))^2(Z^6f_7(T/Z, X/Z))\\
&& = 4Z^3(-4XZ^4r^2 + 4T^3Z^2r - 4T^2Z^3r + 2T^5 + T^4X - 4T^4Z - 6T^3XZ + 3T^2XZ^2\\
&& + 2T^2Z^3 + 2TXZ^3 + 4X^2Z^3 + XZ^4)^7.
\end{eqnarray*}
Since we can choose $r$ arbitrary, $\mcB$ given by $-4Z^9f_7(T/Z, X/Z)= 0$ satisfies 
infinitely many quasi-toric relations of type (2, 7, 2) such that 
\begin{eqnarray*}
h_1&=&Z^{19}b_0(T/Z, X/Z),\\
h_2&=&-4XZ^4r^2 + 4T^3Z^2r - 4T^2Z^3r + 2T^5 + T^4X - 4T^4Z - 6T^3XZ + 3T^2XZ^2\\
&& + 2T^2Z^3 + 2TXZ^3 + 4X^2Z^3 + XZ^4, \\
h_3&=&Z^{16}b_1(T/Z, X/Z), F_1=1, F_2=-4Z^3\,  \mbox{and} \, F_3 =Z^6f_7(T/Z, X/Z).
\end{eqnarray*}
 }
 \end{ex}

\noindent Ai TAKAHASHI and Hiro-o TOKUNAGA\\
Department of Mathematical  Sciences, Graduate School of Science, \\
Tokyo Metropolitan University, 1-1 Minami-Ohsawa, Hachiohji 192-0397 JAPAN \\
{\tt tokunaga@tmu.ac.jp}

\end{document}